\pdfoutput=1
\documentclass[toc]{mathprint}
\usepackage{rutarmacros}
\usepackage{macros}
\title[Overlapping Self-affine Sets]{Assouad-type Dimensions of Overlapping Self-affine Sets}

\author[Fraser]
  {Jonathan M. Fraser}
  {Mathematical Institute, University of St Andrews, St Andrews KY16 9SS, Scotland}
  {jmf32@st-andrews.ac.uk}

\author[Rutar]
  {Alex Rutar}
  {Mathematical Institute, University of St Andrews, St Andrews KY16 9SS, Scotland}
  {alex@rutar.org}

\begin{document}
\begin{abstract}
    We study the Assouad and quasi-Assoaud dimensions of dominated rectangular self-affine sets in the plane.
    In contrast to previous work on the dimension theory of self-affine sets, we assume that the sets satisfy certain separation conditions on the projection to the principal axis, but otherwise have arbitrary overlaps in the plane.
    We introduce and study regularity properties of a certain symbolic non-autonomous iterated function system corresponding to ``symbolic slices'' of the self-affine set.
    We then establish dimensional formulas for the self-affine sets in terms of the dimension of the projection along with the maximal dimension of slices orthogonal to the projection.
    Our results are new even in the case when the self-affine set satisfies the strong separation condition: in fact, as an application, we show that self-affine sets satisfying the strong separation condition can have distinct Assouad and quasi-Assouad dimensions, answering a question of the first named author.
\end{abstract}

\section{Introduction}
Dimension theory is generally concerned with the scaling properties of subsets of a metric space.
The Assouad dimension is a coarse measurement of scaling, in that it captures the worst case scaling behaviour across all locations in the set and across all pairs of scales.
This is in contrast to the Hausdorff dimension, which captures the average scaling of a set at arbitrarily small scales.
We also study the quasi-Assouad dimension, which forces an exponential separation between the scales considered by the Assouad dimension and therefore lies between the Hausdorff and Assouad dimensions.

For highly homogeneous sets---such as Ahlfors regular sets---the Hausdorff, quasi-Assouad, and Assouad dimensions coincide.
However, for many dynamically invariant sets, these dimensions can be different.
There are two particularly notable examples.
Firstly, self-similar sets in $\R$ which satisfy the exponential separation condition (see \cref{d:esc}) but not the open set condition have equal Hausdorff and quasi-Assouad dimension, but Assouad dimension 1 \cite{fy2018a}.
(The authors believe that \emph{all} self-similar sets in $\R$ have equal Hausdorff and quasi-Assouad dimension.)
Secondly, for dominated and irreducible planar self-affine sets satisfying the strong separation condition with Assouad dimension greater than or equal to $1$, it can happen that the Hausdorff dimension and Assouad dimension are distinct \cite[Example~3.3]{bkypreprint}.
This is different than the case where the Assouad dimension is less than $1$, where the corresponding self-affine set is Ahlfors regular \cite[Theorem~1.4]{bkypreprint}.

Much work on the dimension theory of planar self-affine sets assumes either irreducibility and a diophantine property in the plane (for some notable examples, see \cite{bhr2019,bkypreprint,hr2021}), or allows failure of irreducibility but imposes robust separation conditions (such as the open set condition) on the self-affine set as well as its projections \cite{bar2007,fra2014,lg1992,mac2011}.
In this paper, we break both of these assumptions and study systems which do not satisfy an irreducibility hypothesis and also permit large overlaps in the plane.
We will introduce a class of rectangular self-affine sets which satisfy a separation condition which permits certain large overlaps in the projection, and which otherwise permits arbitrary overlaps in the plane.

Our main technique to study Assouad-type dimensions of such systems is through a non-autonomous iterated function system which symbolically encodes the vertical slices of the self-affine set.
The construction and results relating to the ``symbolic slices'' are the key technical contribution of this paper, and we believe that they may be of independent interest.
We prove regularity results for the Assouad-type dimensions of these systems which, when combined with separation in the projection, allows us to establish dimensional results in the plane.
Our results are new even in the case when the self-affine set satisfies separation conditions in the plane.

We were originally motivated to study this class of examples by \cite[Question~17.5.4]{fra2021}, which asks whether or not self-affine sets satisfying the strong separation condition can have distinct quasi-Assouad and Assouad dimensions.
We answer this question by showing that such a phenomenon is possible.

\subsection{Assouad-type dimensions}
Throughout, we work in $\R^d$ equipped with the infinity norm; often, $d=2$.
The choice of infinity norm is merely for convenience, since all other equivalent norms give the same values for dimensions.
Given a bounded set $F\subset\R^d$, the \defn{Assouad dimension} of $F$ is the number
\begin{align*}
    \dimA F=\inf\Bigl\{s\geq 0:(\exists C>0)&(\forall 0 < r\leq R\leq 1)(\forall x\in F)\\
                                           &N_{r}(F\cap B(x,R))\leq C\Bigl(\frac{R}{r}\Bigr)^s\Bigr\}
\end{align*}
where $N_r(E)$ is the least number of balls of radius $r$ required to cover the set $E$.
The Assouad dimension was brought to the forefront in \cite{ass1977} in the context of embedding theory, and has since received a large amount of attention in conformal geometry and fractal geometry.
We refer the readers to the books \cite{fra2021,mt2010,rob2010} for more background and details.

Closely related to the Assouad dimension is the \defn{quasi-Assouad dimension}, introduced by Lü and Xi \cite{lx2016}, and the \defn{Assouad spectrum}, introduced by Fraser and Yu \cite{fy2018b}.
These notions of dimension impose an exponential gap between the upper scale $R$ and the lower scale $r$.
Given a set $F\subset\R^d$ and parameter $\theta\in(0,1)$, the \defn{Assouad spectrum} of $F$ at $\theta$ is
\begin{align*}
    \dimAs\theta F=\inf\Bigl\{s:(\exists C>0)&(\forall 0 < \delta \leq 1)(\forall x\in F)\\
                                                  &N_{\delta^{1/\theta}}(F\cap B(x,\delta))\leq C\Bigl(\frac{\delta}{\delta^{1/\theta}}\Bigr)^s\Bigr\}.
\end{align*}
It always holds that $\dimAs\theta F$ is a uniformly continuous function of $\theta$, so we can define the \defn{quasi-Assouad dimension} by $\dimqA F\coloneqq\lim_{\theta\to 1}\dimAs\theta F$.
Note that this is not the original definition from \cite{lx2016}, but is equivalent by \cite[Theorem~2.1]{fhh+2019}.
Moreover, a wide range of behaviour is possible for the Assouad spectrum.
A discussion of the general properties of the Assouad spectrum, and in particular a complete characterization of possible forms of the Assouad spectrum, is given in \cite{rutpreprint}.

The Assouad spectrum and quasi-Assouad dimension are useful to obtain finer information about the scaling structure of a set in the context of Assouad dimension.
In this note, we study the Assouad and quasi-Assouad dimensions of rectangular dominated planar self-affine sets.
We impose a mild separation condition in the projection (which allows exact overlaps), but notably we require no other assumptions on the matrices and allow any overlaps in the plane.
This class of self-affine sets, along with the relevant assumptions, are defined and discussed in the following section.

\subsection{Iterated function systems and separation conditions}\label{ss:not}
Let $\mathcal{I}$ be a finite index set and consider contraction ratios $\{(\alpha_i,\beta_i)\}_{i\in\mathcal{I}}$ where $0<\beta_i<\alpha_i<1$ for all $i\in\mathcal{I}$.
Fix translations $u_i,v_i\in\R$ for each $i\in\mathcal{I}$.
Then for each $i\in\mathcal{I}$, define the map $T_i\colon\R^2\to\R^2$ by $T_i(x,y)=(\alpha_i x+u_i,\beta_i y+v_i)$.
We refer to the family $\{T_i\}_{i\in\mathcal{I}}$ as a \defn{dominated rectangular self-affine iterated function system}.
Without loss of generality, we may assume that $T_i([0,1]^2)\subset[0,1]^2$.

As usual, let $\mathcal{I}^*=\bigcup_{n=0}^\infty\mathcal{I}^n$ and for $\sigma=(i_1,\ldots,i_n)\in\mathcal{I}^*$, write
\begin{align*}
    T_\sigma&=T_{i_1}\circ\cdots\circ T_{i_n}\\
    \alpha_\sigma &= \alpha_{i_1}\cdots \alpha_{i_n}\\
    \beta_\sigma &= \beta_{i_1}\cdots \beta_{i_n}.
\end{align*}
We also write $\beta_{\max}=\max\{\beta_i:i\in\mathcal{I}\}$, and similarly define $\beta_{\min}$, $\alpha_{\max}$, and $\alpha_{\min}$.

Let $\pi\colon\R^2\to\R$ denote the projection $\pi(x,y)=x$.
For $\sigma\in\mathcal{I}^n$, we let $S_\sigma\colon\R\to\R$ denote the unique similarity satisfying $S_\sigma\circ\pi=\pi\circ T_\sigma$.
We then define the \defn{projected IFS} by $\{S_i\}_{i\in\mathcal{I}}$.
Observe that $S_i(x)=\alpha_i x+u_i$ for $i\in\mathcal{I}$.

We will impose some separation conditions on the IFS $\{T_i\}_{i\in\mathcal{I}}$ by constraining the projected IFS $\{S_i\}_{i\in\mathcal{I}}$.
    Fix an IFS $\{S_i\}_{i\in\mathcal{I}}$ of similarities with attractor $K$, i.e.~$S_i(x)=\alpha_i x+u_i$ for $\alpha_i\in(0,1)$.
\begin{definition}\label{d:wsc-defs}
    For $r\in(0,1)$, let
    \begin{equation*}
        \Lambda_r\coloneqq\{(i_1,\ldots,i_n)\in\mathcal{I}^*:\alpha_{i_1}\cdots \alpha_{i_n}\leq r<\alpha_{i_1}\cdots \alpha_{i_{n-1}}\}
    \end{equation*}
    and
    \begin{equation*}
        t_r\coloneqq \sup_{x\in K}\#\bigl\{S_\sigma:\sigma\in\Lambda_r,S_\sigma(K)\cap B(x,r)\neq\emptyset\bigr\}.
    \end{equation*}
    We say that the IFS $\{S_i\}_{i\in\mathcal{I}}$ satisfies the \defn{weak separation condition} (WSC) if $(t_r)_{r>0}$ is bounded, and the \defn{asymptotically weak separation condition} (AWSC) if $\lim_{r\to 0}\frac{\log t_r}{\log r}=0$.
    Finally, we say that the IFS satisfies the \defn{open set condition} (OSC) if it satisfies the WSC and $S_\sigma\neq S_\tau$ for all $n\in\N$ and $\sigma\neq\tau\in\mathcal{I}^n$
\end{definition}
The above definition of the OSC is equivalent to the usual notion with respect to an open set by \cite{sch1994}.

Next, we define the \defn{exponential separation condition}, which was introduced and studied in \cite{hoc2014}.
\begin{definition}\label{d:esc}
    Fix an IFS of similarities $\{S_i\}_{i\in\mathcal{I}}$.
    Given words $\sigma,\tau\in\mathcal{I}^n$, let
    \begin{equation*}
        d(\sigma,\tau)=
        \begin{cases}
            |S_\sigma(0)-S_\tau(0)| &: \alpha_\sigma=\alpha_\tau\\
            \infty &:\text{ otherwise}.
        \end{cases}
    \end{equation*}
    Let $\Delta_n=\min_{\sigma\neq\tau\in\mathcal{I}^n}d(\sigma,\tau)$.
    We say that the IFS $\{S_i\}_{i\in\mathcal{I}}$ satisfies the \defn{exponential separation condition} (ESC) if $\liminf_{n\to\infty}(-\log\Delta_n)/n<\infty$.
\end{definition}
The AWSC holds very generally, as is specified by the following result.
This is \cite[Theorem~1.3]{bf2021}, but follows essentially from Shmerkin's result on $L^q$-dimensions of self-similar sets \cite{shm2019}.
\begin{proposition}[\cite{bf2021}]\label{p:esc-class}
    Let $\{S_i\}_{i\in\mathcal{I}}$ be an IFS of similarities in $\R$ with contraction ratios $\{r_i\}_{i\in\mathcal{I}}$.
    Let $s$ denote the unique solution to $\sum_{i\in\mathcal{I}}r_i^s=1$.
    Suppose $\{S_i\}_{i\in\mathcal{I}}$ satisfies the ESC.
    Then $\{S_i\}_{i\in\mathcal{I}}$ satisfies the AWSC if and only if $s\leq 1$.
\end{proposition}
In particular, within many parametrized families of IFSs, if the similarity dimension in the parametrized family is uniformly bounded above by $1$, then the AWSC holds outside of a small exceptional set (typically of Hausdorff dimension 0)---see, for example, \cite[Theorem~1.8]{hoc2014}.
In general, it is not known if the ESC assumption can be replaced with the assumption ``has no exact overlaps''.

Unlike the ESC, however, the AWSC also permits exact overlaps.
For example, any IFS satisfying the WSC satisfies the AWSC.
The AWSC is also known to hold for explicit examples of IFSs with exact overlaps which do not satisfy the WSC, such as Bernoulli convolutions with Salem numbers.
As another example, it is known that the IFS $\{x\mapsto x/3,x\mapsto x/3+t,x\mapsto x/3+2/3\}$ satisfies the AWSC for all $t\in(0,1)$ by \cite[Theorem~1.5]{bf2021} and \cite[Theorem~1.6]{bf2021} with \cref{p:esc-class}.
However, the WSC fails when $t$ is irrational by \cite[Theorem~2]{ken1997}, since any self-similar set $K$ satisfying the WSC has positive $\dimH K$-dimensional Hausdorff measure \cite[Theorem~2]{fal1989}.

\subsection{Symbolic fibre dimensions}\label{ss:s-fibre}
Let $\mathcal{S}$ denote the semigroup under composition with generators $\{S_i\}_{i\in\mathcal{I}}$.
For $f\in\mathcal{S}$, let
\begin{equation*}
    \mathcal{J}_f=\{\sigma\in\mathcal{I}^*:S_\sigma=f\}
\end{equation*}
and we define the \defn{$f$-fibred IFS} $\{F_\sigma\}_{\sigma\in\mathcal{J}_f}$ where $F_\sigma\circ\overline{\pi}=\overline{\pi}\circ T_\sigma$ and $\overline{\pi}\colon\R^2\to\R$ is the projection $\overline{\pi}(x,y)=y$.
Note that $\{F_\sigma\}_{\sigma\in\mathcal{J}_f}$ are IFSs of similarities in $\R$ for all $f\in\mathcal{S}$.

Let $\Omega$ denote the collection of sequences $\eta=(\eta_n)_{n=1}^\infty\subset\mathcal{S}$ where for each $n\in\N$ there is some $i\in\mathcal{I}$ so that $\eta_{n}=\eta_{n-1}\circ S_i$ (writing $\eta_0=\id$).
If $\mathcal{S}$ is a free semigroup, then $\mathcal{I}^{\N}\cong\Omega$ by the map $(i_n)_{n=1}^\infty\mapsto (S_{i_1},S_{i_1}\circ S_{i_2},\ldots)$.
Let $X\subset\R$ be a non-empty compact set and for $\eta\in\Omega$, set
\begin{equation}\label{e:fibre-def}
    E_{n,\eta}\coloneqq\bigcup_{\sigma\in\mathcal{J}_{\eta_n}}F_\sigma(X)
\end{equation}
and let $E_\eta=\lim_{n\to\infty}E_{n,\eta}$ with respect to the Hausdorff metric on compact sets.
We will see that the limit exists and does not depend on the choice of $X$ in \cref{l:symb-fibre-limit}.
Observe that if $x_\eta=\lim_{n\to\infty}\eta_n(0)$, then $\{x_\eta\}\times E_\eta\subset K$.
In this sense, one can interpret the space $\Omega$ as a symbolic analogue of $\pi(K)$ and the set $E_\eta$ as a symbolic vertical slice of $K$.

Our main goal is to prove dimension formulas for $\dimA K$ and $\dimAs\theta K$ (for $\theta\in(0,1)$ sufficiently close to 1) in terms of the dimension of the projection and a \defn{symbolic fibre dimension}
\begin{equation*}
    \max_{\eta\in\Omega}\dim E_\eta
\end{equation*}
where $\dim$ denotes either the Assouad dimension or the Assouad spectrum at some $\theta$ (that the maximum is attained is given in \cref{p:uniform-A}).

Note that if $\mathcal{S}$ is a free semigroup and the $S_i$-fibred IFSs satisfy the OSC with respect to $(0,1)$ for all $i\in\mathcal{I}$, then $E_\eta$ is the attractor of a non-autonomous self-similar IFS, as introduced and studied in \cite{ru2015}.
In this situation, it holds that
\begin{equation*}
    \max_{\eta\in\Omega}\dimA E_\eta=\max\{a_i\}_{i\in\mathcal{I}}
\end{equation*}
where $a_i$ is the Hausdorff dimension of the attractor of the $S_i$-fibred IFS for $i\in\mathcal{I}$.
For instance, this is the case for Gatzouras--Lalley carpets (their Assouad dimensions were originally computed in \cite{mac2011}).
However, we emphasize that we make no assumptions on the overlaps within the fibres.
Even under the WSC, subtle behaviour is possible.
For $0<\lambda\leq (3-\sqrt{5})/2\approx 0.382$, let
\begin{equation*}
    \Phi_\lambda = \{x\mapsto \lambda x,x\mapsto \lambda x+\lambda-\lambda^2,x\mapsto \lambda x+(1-\lambda)\}.
\end{equation*}
Note that $\Phi_\lambda$ satisfies the WSC \cite[Example~2.3]{lnr2001} and its attractor has dimension $(\log\frac{3+\sqrt{5}}{2})/(-\log\lambda)$ \cite[Example~5.4]{nw2001} for all such $\lambda$.
On the other hand, one can choose $\lambda_1$ and $\lambda_2$ arbitrarily small so that the IFS
\begin{equation}\label{e:lambda-ifs}
    \Phi=\{f\circ g:f\in\Phi_{\lambda_1},g\in\Phi_{\lambda_2}\}
\end{equation}
does not satisfy the WSC (see \cref{r:it-wsc-fail} for the details).
In particular, for $\lambda_1,\lambda_2$ chosen as above, the IFS with maps
\begin{align*}
    T_1(x,y)&=\bigl(\frac{x}{2}, \lambda_1 y+(1-\lambda_1)\bigr) & T_2(x,y)&=\bigl(\frac{x}{2}+\frac{1}{2}, \lambda_2 y+(1-\lambda_2)\bigr)\\
    T_3(x,y)&=\bigl(\frac{x}{2}, \lambda_1 y+\lambda_1-\lambda_1^2\bigr) & T_4(x,y)&=\bigl(\frac{x}{2}+\frac{1}{2}, \lambda_2 y+\lambda_2-\lambda_2^2\bigr)\\
    T_5(x,y)&=\bigl(\frac{x}{2}, \lambda_1 y\bigr) & T_6(x,y)&=\bigl(\frac{x}{2}+\frac{1}{2}, \lambda_2 y\bigr)
\end{align*}
has $\max_{\eta\in\Omega}\dimA E_\eta=1$ so, by \cref{it:main-A} below, the Assouad dimension of the attractor is $2$.
On the other hand, if $\lambda_1=\lambda_2=\lambda$, then the attractor is a product of self-similar sets and has Assouad dimension $1+(\log\frac{3+\sqrt{5}}{2})/(-\log\lambda)$ by \cref{c:ss-prod}.

\subsection{Statement and discussion of results}
Our first main contribution is the following formula for the Assouad dimension for a class of IFSs of affinities in terms of the dimensions of the principal projection and the maximal symbolic dimension of a fibre.
We can also transfer the symbolic dimension of a fibre to the dimension of a vertical slice.
\begin{itheorem}\label{it:main-A}
    Let $\{T_i\}_{i\in\mathcal{I}}$ be a dominated rectangular IFS of affinities with attractor $K$.
    Then
    \begin{equation}\label{e:dimA-formula-intro}
        \dimA K\geq \dimA \pi(K)+\max_{\eta\in\Omega}\bigl(\dimA E_\eta\bigr).
    \end{equation}
    Moreover, if the projected IFS $\{S_i\}_{i\in\mathcal{I}}$ satisfies the WSC, then equality holds in \cref{e:dimA-formula-intro} and
    \begin{equation}\label{e:fibre-form-intro}
        \dimA K=\dimA \pi(K) + \max_{x\in\pi(K)}\dimA\bigl(\pi^{-1}(x)\cap K\bigr).
    \end{equation}
\end{itheorem}
The proof of this result is given in \cref{l:fibre} and \cref{t:dima}, and the fact that the supremum in \cref{e:dimA-formula-intro} is always attained as a maximum is given in \cref{p:uniform-A}.
The lower bound for the Assouad dimension is given by constructing a particular weak tangent.
To do this, we observe a simple but important tangent regularity result for self-similar sets (this is given in \cref{l:wt-choice}).
As a quick application of this, we also show that products of compact sets with self-similar sets have maximal Assouad dimension: see \cref{c:ss-prod}.

This also establishes the fibre formula discussed in \cite[Question~17.5.1]{fra2021}, for the systems considered in the theorem.
We believe this formula to hold more generally for any genuinely self-affine set, but this general problem seems to be beyond our current techniques.

We also note the relationship of these results with \cite{fj2017}, though that paper assumes the OSC in the plane, regularity results for the projection to the line, and the homogeneity assumption that the contraction ratios satisfy $\alpha_i=\alpha$ and $\beta_i=\beta$ for some constants $\alpha,\beta$ and all $i\in\mathcal{I}$.
There, the local dimension plays a similar role to the symbolic fibre dimension in \cref{it:main-A}.
In certain cases, their results imply the formula given in \cref{e:fibre-form-intro}.

Of course, it always holds that $\dimA K\leq\dimA \pi(K)+1$ since $K\subseteq\pi(K)\times[0,1]$.
In particular, if $\dimA E_{\eta}$ can be arbitrarily close to $1$, then equality also holds in \cref{e:dimA-formula-intro}.
Moreover, this implies that equality holds in \cref{e:fibre-form-intro} as well.
This gives the following immediate application:
\begin{icorollary}\label{ic:Assouad-formula}
    Let $\{T_i\}_{i\in\mathcal{I}}$ be a dominated rectangular IFS of affinities with attractor $K$.
    Suppose $\max_{\eta\in\Omega}\bigl(\dimA E_\eta\bigr)=1$.
    Then
    \begin{equation*}
        \dimA K=\dimA\pi(K)+1=\dimA\pi(K)+\max_{x\in\pi(K)}\dimA\bigl(\pi^{-1}(x)\cap K\bigr).
    \end{equation*}
\end{icorollary}
For dominated rectangular IFSs of affinities in general, one should not expect a symbolic formula such as \cref{e:dimA-formula-intro} to hold: see \cref{r:awsc-req} for more details.

It is also proven in \cite[Theorem~1.3]{fhor2015} that for an IFS of similarities in $\R$ with attractor $F$, if $\dimA F<1$, then the IFS satisfies the WSC.
Moreover, since $\alpha_i>\beta_i$ for some $i\in\mathcal{I}$, by considering compositions of the form $S_i\circ\cdots\circ S_i(K)$, we observe that $\pi(K)$ is always a weak pseudo-tangent for $K$ and so that $\dimA K\geq\dimA\pi(K)$ (see \cref{ss:tan}).
(In general, the Assouad dimension of rectangular self-similar sets in the plane can increase under projection \cite[Section~3.1]{fra2014}.)
In particular, if $\dimA K<1$, then $\pi(K)$ is the attractor of an IFS satisfying the WSC, which provides the following application:
\begin{icorollary}
    Let $\{T_i\}_{i\in\mathcal{I}}$ be a dominated rectangular IFS of affinities with attractor $K$.
    Suppose $\dimA K<1$.
    Then
    \begin{equation*}
        \dimA K=\dimA \pi(K) + \max_{x\in\pi(K)}\dimA\bigl(\pi^{-1}(x)\cap K\bigr).
    \end{equation*}
\end{icorollary}

Our second main contribution is a formula for the Assouad spectrum for sufficiently large $\theta$ in terms of the principal projection and the maximal symbolic dimension of a fibre.
\begin{itheorem}\label{it:main-As}
    Let $\{T_i\}_{i\in\mathcal{I}}$ be a dominated rectangular IFS of affinities with attractor $K$.
    Then for all $\max_{i\in\mathcal{I}}\frac{\log\alpha_i}{\log\beta_i}\leq\theta<1$,
    \begin{equation}\label{e:dimAs-formula-intro}
        \dimAs\theta K\geq\dimH \pi(K)+\max_{\eta\in\Omega}\bigl(\dimAs\theta E_\eta\bigr)
    \end{equation}
    Moreover, if the projected IFS $\{S_i\}_{i\in\mathcal{I}}$ satisfies the AWSC, then equality holds in \cref{e:dimAs-formula-intro}.
\end{itheorem}
This proof is given in \cref{t:dimas}.
The exponential separation of scales in the Assouad spectrum allows us to ignore sub-exponential errors from the AWSC.
This is in contrast to the Assouad dimension, which is affected by sub-exponential growth of constants: this is essentially a manifestation in the plane of the dichotomy result for Assouad dimensions of self-similar sets \cite[Theorem~1.3]{fhor2015}.

Finally, since it is well-known that sets with arbitrarily small similarity dimension satisfying the AWSC but not the WSC exist, we obtain the following corollary:
\begin{icorollary}
    For every $\epsilon>0$, there is an IFS of affinities $\{T_i\}_{i\in\mathcal{I}}$ satisfying the strong separation condition with attractor $K$ such that $\dimqA K\leq\epsilon$ and $\dimA K\geq 1$.
\end{icorollary}
See \cref{ss:ques} for an explicit construction accompanied by a detailed discussion.
In particular, this answers \cite[Question~17.5.4]{fra2021}.

\section{Proofs of main results}
\subsection{Tangents and product sets}\label{ss:tan}
Given a set $E\subset\R^d$, denote the $\delta$-neighbourhood of $E$ by
\begin{equation*}
    E^{(\delta)}=\{x\in\R^d:\norm{x-y}<\delta\text{ for some }y\in E\}.
\end{equation*}
Given non-empty compact sets $E$ and $F$, define the \defn{Hausdorff pseudo-distance} by
\begin{equation*}
    p_{\mathcal{H}}(E,F)=\inf\{\delta>0:E\subset F^{(\delta)}\}
\end{equation*}
and the \defn{Hausdorff distance}
\begin{equation*}
    d_{\mathcal{H}}(E,F)=\max\{p_{\mathcal{H}}(E,F),p_{\mathcal{H}}(F,E)\}.
\end{equation*}
If $X\subset\R^d$ is a compact set, the set $\mathcal{K}(X)$ of all non-empty compact subsets of $X$ equipped with the Hausdorff distance $d_{\mathcal{H}}$ is a compact metric space.

Given a similarity $W\colon\R^d\to\R^d$, we may write
\begin{equation*}
    W(\bm{x})=\gamma O\bm{x}+\bm{a}
\end{equation*}
where $O\in\mathcal{O}(\R^d)$ is an orthogonal matrix, $\gamma>0$ is a constant, and $\bm{a}\in\R^d$.
We refer to $\gamma$ as the \defn{scaling ratio} of $W$.

The notion of a weak tangent was introduced in \cite{mt2010}, with ideas going back to \cite{kl2004}.
We also find it convenient to use the slightly modified notion of a weak pseudo-tangent, which was introduced in \cite{fhor2015}.
Let $F$ and $\widehat{F}$ be compact subsets of $\R^d$.
We say that $\widehat{F}$ is a \defn{weak pseudo-tangent} of $F$ if there exists a sequence of similarities $(T_k)_{k=1}^\infty$ with scaling ratios diverging to infinity such that $\lim_{k\to\infty}p_{\mathcal{H}}(\widehat{F},T_k(F))=0$.
Similarly, we say that $\widehat{F}$ is a \defn{weak tangent} of $F$ if $\widehat{F}$ is a weak pseudo-tangent and additionally $\lim_{k\to\infty}T_k(F)\cap B(0,1)=\widehat{F}$ in the Hausdorff distance.

We recall the following result, which is \cite[Proposition~3.7]{fhor2015} and also essentially follows from \cite[Proposition~6.1.5]{mt2010}.
\begin{lemma}[\cite{fhor2015,mt2010}]\label{l:tan-lower}
    If $\widehat{F}$ is a weak pseudo-tangent of $F$, then $\dimA\widehat{F}\leq\dimA F$.
\end{lemma}
We also recall the following strong converse, which is given in \cite[Proposition~5.7]{kor2018} (see also \cite[Theorem~5.1.3]{fra2021}).
\begin{proposition}[\cite{kor2018}]\label{p:wt-exist}
    Let $F\subset\R^d$ be closed and non-empty with $\dimA F=s$.
    Then there is a compact set $E\subset\R^d$ such that $\mathcal{H}^s(E)>0$ and $E$ is a weak tangent of $F$.
\end{proposition}
We begin with the following simple, but useful, observation concerning weak tangents of self-similar sets, which was made by S.~Troscheit:
\begin{lemma}\label{l:wt-choice}
    Let $\{S_i\}_{i\in\mathcal{I}}$ be an IFS of similarities with attractor $K$ and let $s=\dimA K$.
    Then $K$ has a weak pseudo-tangent $\widehat{K}$ with $\mathcal{H}^s(\widehat{K})>0$ with respect to a sequence of similarities $(U_j)_{j=1}^\infty$ with scaling ratios $\gamma_j$ satisfying $\gamma_j\leq\gamma_{j+1}\leq C\gamma_j$ for some constant $C>0$.
\end{lemma}
\begin{proof}
    By \cref{p:wt-exist}, with $s=\dimA K$, there is a compact $\widehat{K}$ with $\mathcal{H}^s(\widehat{K})>0$ and a sequence of similarities $(\widehat{U}_j)_{j=1}^\infty$ with scaling ratios $\widehat{\gamma}_j$ such that
    \begin{equation*}
        \lim_{j\to\infty}d_{\mathcal{H}}(\widehat{K},\widehat{U}_j(K)\cap B(0,1))=0
    \end{equation*}
    and $(\widehat{\gamma}_j)_{j=1}^\infty$ converges monotonically to infinity.

    Fix some index $i_0\in\mathcal{I}$ and consider the similarity $\widehat{U}_j\circ S_{i_0}^{-1}$, which has scaling ratio $\widehat{\gamma}_j\cdot r_{i_0}^{-1}$.
    Moreover, since $S_{i_0}(K)\subset K$,
    \begin{equation*}
        \widehat{U}_j\circ S_{i_0}^{-1}(K)\supset \widehat{U}_j(K)\cap B(0,1).
    \end{equation*}
    In particular, $p_{\mathcal{H}}(\widehat{K},\widehat{U}_j\circ S_{i_0}^{-1}(K))\leq d_{\mathcal{H}}(\widehat{K},\widehat{U}_j(K)\cap B(0,1))$.
    Now for each $j\in\N$, let $m_j$ be maximal such that $\widehat{\gamma}_j r_{i_0}^{-m_j}<\widehat{\gamma}_{j+1}$.
    Then $\widehat{K}$ is a weak pseudo-tangent of $K$ with respect to the sequence of similarities
    \begin{equation*}
        (\widehat{U}_1,\widehat{U}_1\circ S_{i_0}^{-1},\ldots,\widehat{U}_1\circ S_{i_0}^{-m_1},\widehat{U}_2,\widehat{U}_2\circ S_{i_0}^{-1},\ldots),
    \end{equation*}
    which satisfy the required properties.
\end{proof}
We note the following quick application:
\begin{corollary}\label{c:ss-prod}
    Let $K\subset\R^d$ be a self-similar set and let $F\subset\R^\ell$ be non-empty and compact.
    Then $\dimA K\times F=\dimA K+\dimA F$.
\end{corollary}
\begin{proof}
    It is well-known (and easy to show) that $\dimA K\times F\leq\dimA K+\dimA F$.
    To obtain the lower bound, we construct an appropriate weak pseudo-tangent.
    Applying \cref{p:wt-exist}, get a weak tangent $\widehat{F}$ with $\mathcal{H}^{\dimA F}(\widehat{F})>0$ for $F$ with respect to a sequence similarities $(U_j)_{j=1}^\infty$.
    Let
    \begin{equation*}
        U_j(\bm{x})=\gamma_j O_j\bm{x}+\bm{u}_j
    \end{equation*}
    for orthogonal matrices $O_j$.
    By \cref{l:wt-choice}, there is some $C>0$ such that $K$ has a weak pseudo-tangent $\widehat{K}$ with $\mathcal{H}^{\dimA K}(\widehat{K})>0$ with respect to a sequence of similarities $(\Phi_j)_{j=1}^\infty$ where
    \begin{equation*}
        \Phi_j(\bm{x})=\eta_j Q_j\bm{x}+\bm{v}_j
    \end{equation*}
    for orthogonal matrices $Q_j$ and $\eta_{j}\leq \eta_{j+1}\leq C\eta_j$.
    In particular, by passing to a subsequence if necessary, we may assume that $\eta_j\in[\gamma_j, C\gamma_j]$.

    Then with $\Psi_j:\R^{d+\ell}\to\R^{d+\ell}$ given by
    \begin{equation*}
        \Psi_j(\bm{x},\bm{y})=\gamma_j\begin{pmatrix}Q_j&0\\0&O_j\end{pmatrix}\begin{pmatrix}\bm{x}\\\bm{y}\end{pmatrix}+\begin{pmatrix}\bm{u}_j\\\bm{v}_j\end{pmatrix},
    \end{equation*}
    we observe that there are similarities $f_j$ and $g_j$ so that
    \begin{equation*}
        \lim_{j\to\infty}p_{\mathcal{H}}\bigl(f_j(\widehat{K})\times g_j(\widehat{F}),\Psi_j(K\times F)\bigr)=0
    \end{equation*}
    where $f_j$ and $g_j$ have scaling ratios in the interval $[C^{-1},C]$.
    Thus passing to a subsequence and using compactness of the groups of orthogonal matrices $\mathcal{O}(\R^{d})$ and $\mathcal{O}(\R^\ell)$, we conclude that there are similarities $f$ and $g$ so that $f(\widehat{K})\times g(\widehat{F})$ is a weak pseudo-tangent for $K\times F$.
    Moreover,
    \begin{equation*}
        \dimH f(\widehat{K})\times g(\widehat{F})\geq \dimA K+\dimA F.
    \end{equation*}
    Thus $\dimA K\times F\geq\dimA K+\dimA F$.
\end{proof}
\begin{remark}
    In general, $\dimA K\times F\geq\max\{\dimA K,\dimA F\}$ since the Assouad dimension is monotonic under inclusion.
    This is sharp: using a homogeneous Moran set construction in which the scales at which the sets $K$ and $F$ are large are complementary, one can construct compact sets $K$ and $F$ with arbitrary Assouad dimension so that $\dimA K\times F=\max\{\dimA K,\dimA F\}$.
    Such constructions are described in \cite{ors2016}.
\end{remark}

\subsection{Symbolic fibres and uniformity of Assouad-type dimensions}
Throughout this section, we fix a dominated rectangular IFS of affinities $\{T_i\}_{i\in\mathcal{I}}$.
We recall that the notation used in this section is established in \cref{ss:s-fibre}.

We first prove existence of the limits defining the symbolic fibres.
For $k\in\N$ and $\eta=(\eta_n)_{n=1}^\infty\in\Omega$, we write $\eta|k=\eta_k\circ \eta_{k-1}^{-1}\in\mathcal{S}$ (taking $\eta_0=\id$).
\begin{lemma}\label{l:symb-fibre-limit}
    Let $X\subset\R$ be compact and non-empty.
    Then for all $\eta\in\Omega$, writing
    \begin{equation}\label{e:E-step-def}
        E_{k,\eta}=\bigcup_{\sigma\in\mathcal{J}_{\eta_k}}F_\sigma(X),
    \end{equation}
    the limit $E_\eta=\lim_{k\to\infty}E_{k,\eta}$ exists and is independent of the choice of $X$.
    Moreover, there is a constant $C>0$ so that
    \begin{equation}\label{e:ph-dist}
        p_{\mathcal{H}}(E_{k,\eta},E_\eta)\leq C\max\{\beta_\sigma:\sigma\in\mathcal{J}_{\eta_k}\}.
    \end{equation}
\end{lemma}
\begin{proof}
    Since the maps $F_i$ for $i\in\mathcal{I}$ are strictly contracting similarities, it is clear that the convergence does not depend on the choice of $X$.
    Thus without loss of generality, we may take $X$ compact so that $F_i(X)\subset X$ for all $i\in\mathcal{I}$.
    Recall that there is a sequence $(S_{i_n})_{n=1}^\infty$ with $i_n\in\mathcal{I}$ for all $n$ such that $\eta=(S_{i_1}\circ\cdots\circ S_{i_n})_{n=1}^\infty$.
    We now define
    \begin{equation*}
        A_{k,n}\coloneqq\bigcup_{\sigma\in\mathcal{J}_{\eta_k}}\bigcup_{(\tau_1,\ldots,\tau_n):\tau_i\in\mathcal{J}_{S_{k+i}}}F_\sigma\circ F_{\tau_1}\circ\cdots\circ F_{\tau_n}(X)\qquad\text{and}\qquad A_k\coloneqq\bigcap_{n=1}^\infty A_{k,n}.
    \end{equation*}
    We make a few observations:
    \begin{itemize}[nl]
        \item Since $F_i(X)\subset X$ for all $i\in\mathcal{I}$, the sequence $(A_{k,n})_{n=1}^\infty$ is a nested sequence of compact sets, so $A_k$ is compact and non-empty.
        \item If $f_1,f_2\in\mathcal{S}$, then for any $\tau_1\in\mathcal{J}_{f_1}$ and $\tau_2\in\mathcal{J}_{f_2}$, $\tau_1\tau_2\in\mathcal{J}_{f_1\circ f_2}$.
            In particular, $A_k\subset A_{k+1}$ for all $k\in\N$.
        \item Since $d_{\mathcal{H}}(E_{k,\eta},A_k)\leq(\diam X)\max\{\beta_\sigma:\sigma\in\mathcal{J}_{\eta_k}\}$, $\lim_{k\to\infty}d_{\mathcal{H}}(E_{k,\eta},A_k)=0$.
    \end{itemize}
    It follows that $\lim_{k\to\infty}A_k=\lim_{k\to\infty}E_{k,\eta}$ exists and is non-empty, and moreover for all $k\leq m$,
    \begin{equation*}
        p_{\mathcal{H}}(E_{k,\eta},E_{m,\eta})\leq (\diam X)\max\{\beta_\sigma:\sigma\in\mathcal{J}_{\eta_k}\}.
    \end{equation*}
    Taking a limit in $m$ yields the claimed \cref{e:ph-dist}.
\end{proof}
We also observe the following ``self-similarity'' of symbolic fibres.
This observation is the key feature which allows us to prove uniform bounds in \cref{p:uniform-A}.
\begin{lemma}\label{l:symb-self-similar}
    Let $\eta\in\Omega$ be arbitrary and let $k\in\N$.
    Then if $\overline{\eta}\in\Omega$ is defined by $\overline{\eta}|n=\eta|(k+n)$ for all $n\in\N$, if $\omega\in\mathcal{J}_{\eta_k}$ is arbitrary,
    \begin{equation*}
        F_\omega(E_{\overline{\eta}})\subset E_\eta.
    \end{equation*}
\end{lemma}
\begin{proof}
    Let $(i_n)_{n=1}^\infty\subset\mathcal{I}$ be chosen so that $\eta_n=S_{i_1}\circ\cdots\circ S_{i_n}$ for all $n\in\N$.
    Note that $\overline{\eta}_n=S_{i_{k+1}}\circ\cdots\circ S_{i_{k+n}}$.
    Fix $\omega\in\mathcal{J}_{\eta_k}$ so $S_\omega=\eta_k$.
Note that $\omega\sigma\in\mathcal{J}_{\eta_{k+n}}$ for any $\sigma\in\mathcal{J}_{\overline{\eta}_k}$ since $S_{\omega}\circ S_{\sigma}=(S_{i_1}\circ\cdots\circ S_{i_k})\circ(S_{i_{k+1}}\circ\cdots\circ S_{i_{k+n}})$.
    Thus if $X\subset\R$ is compact and non-empty, for any $n\in\N$
    \begin{equation*}
        F_\omega(E_{n,\overline{\eta}})=\bigcup_{\sigma\in\mathcal{J}_{\overline{\eta}_k}}F_{\omega\sigma}(X))\subset\bigcup_{\tau\in\mathcal{J}_{\eta_{k+n}}}F_\tau(X)=E_{n+k,\eta}
    \end{equation*}
    where the sets $E_{n,\overline{\eta}}$ and $E_{n,\eta}$ are defined in \cref{e:E-step-def}.
    Passing to the limit in the Hausdorff distance yields the desired result.
\end{proof}
Our main result in this section is the following uniformity result for the symbolic fibre Assouad dimensions.
We note that self-similarity is important: we essentially use the idea from \cref{l:wt-choice} to ``align'' parts of $\Omega$ which have large covering numbers between pairs of scales.
\begin{proposition}\label{p:uniform-A}
    Let
    \begin{equation}\label{e:dimA-def}
        s=\sup_{\eta\in\Omega}\dimA E_\eta.
    \end{equation}
    Then for every $\epsilon>0$, there exists $C_\epsilon>0$ such that for all $0<r\leq R<1$, $\eta\in\Omega$, and $x\in E_\eta$:
    \begin{equation}\label{e:uniform-Assouad}
        N_r(B(x,R)\cap E_\eta)\leq C_\epsilon\Bigl(\frac{R}{r}\Bigr)^{s+\epsilon}.
    \end{equation}
    Moreover, the supremum in \cref{e:dimA-def} is attained as a maximum.
\end{proposition}
\begin{proof}
    Suppose for contradiction that there is some $\epsilon>0$ so that \cref{e:uniform-Assouad} fails, and get sequences $(\eta^{(n)})_{n=1}^\infty\subset\Omega$, $x_n\in E_{\eta^{(n)}}$, $(r_n)_{n=1}^\infty$, $(R_n)_{n=1}^\infty$, and $(C_n)_{n=1}^\infty$ diverging to infinity so that
    \begin{equation}\label{e:constant-failure}
        N_{r_n}(B(x_n,R_n)\cap E_{\eta^{(n)}})\geq C_n\Bigl(\frac{R_n}{r_n}\Bigr)^{s+\epsilon}
    \end{equation}
    for each $n\in\N$.
    For each $n$, let $k_n\in\N$ be sufficiently large so that
    \begin{equation}\label{e:kn-dist}
        d_{\mathcal{H}}(E_{k_n,\eta^{(n)}},E_{\eta^{(n)}})\leq r_n.
    \end{equation}
    For each $n\in\N$, let $g_{n,1},\ldots,g_{n,k_n}$ be the similarity maps in $\{S_i\}_{i\in\mathcal{I}}$ defining $\eta^{(n)}_{k_n}$; that is, for each $1\leq j\leq k_n$, $\eta^{(n)}_{j}=g_{n,1}\circ\cdots\circ g_{n,j}$.
    Then let $\xi\in\Omega$ be the sequence corresponding to the sequence of similarity maps
    \begin{equation*}
        (g_{1,1},\ldots,g_{1,k_1},g_{2,1},\ldots,g_{2,k_2},\ldots).
    \end{equation*}

    Fix $n\in\N$, let $m=k_1+\cdots+k_{n-1}$, and let $\omega\in\mathcal{J}_{\xi_m}$ be arbitrary.
    Observe that if $\overline{\xi}\in\Omega$ satisfies $\overline{\xi}|j=\xi|(j+m)$ for all $j\in\N$, then $F_\omega(E_{\overline\xi})\subset E_\xi$ by \cref{l:symb-self-similar}.
    But now since $k_n$ is chosen to satisfy \cref{e:kn-dist}, we observe that $N_{r_n}(B(\overline x_n,2R_n)\cap E_{\overline \xi}) \geq C_n(R_n/r_n)^{s+\epsilon}$ for some $\overline x_n\in E_{\overline\xi}$ so that
    \begin{equation*}
        N_{\beta_\omega r_n}(B(F_\omega(\overline x_n),\beta_\omega R_n)\cap E_{\xi})\geq \delta C_n\Bigl(\frac{\beta_\omega R_n}{\beta_\omega r_n}\Bigr)^{s+\epsilon}
    \end{equation*}
    for some fixed constant $\delta>0$.
    But $C_n$ diverges to infinity, so $\dimA E_\xi\geq s+\epsilon$, contradicting the definition of $s$.

    To see that the supremum in \cref{e:dimA-def} is attained as a maximum, we simply observe that the same argument applied along a sequence $s-\epsilon_n$ where $\epsilon_n$ converges to zero, rather than the constant sequence $s+\epsilon$, provides some $\xi\in\Omega$ so that $\dimA E_\xi\geq s$.
\end{proof}
In fact, a similar proof gives the analogous result for the Assouad spectrum.
\begin{proposition}\label{p:uniform-qA}
    Let $\theta\in(0,1)$ be arbitrary and let
    \begin{equation}\label{e:dimAs-def}
        s_\theta=\sup_{\eta\in\Omega}\dimAs\theta E_\eta.
    \end{equation}
    Then for every $\epsilon>0$, there exists $C_{\epsilon,\theta}>0$ such that for all $0<R<1$, $\eta\in\Omega$, and $x\in E_\eta$:
    \begin{equation*}
        N_{R^{1/\theta}}(B(x,R)\cap E_\eta)\leq C_{\epsilon,\theta}\Bigl(\frac{R}{R^{1/\theta}}\Bigr)^{s_\theta+\epsilon}.
    \end{equation*}
    Moreover, the supremum in \cref{e:dimAs-def} is attained as a maximum.
\end{proposition}
\begin{proof}
    Follow the proof of \cref{p:uniform-A}, but now the choice of $\eta^{(n)}$, $r_n=R_n^{1/\theta}$, $x_n$, and $C_n$ must be chosen so that, having taken $m=k_1+\cdots+k_{n-1}$ and $\omega\in\mathcal{J}_{\xi_m}$, that $\theta_n$ converges to $\theta$ where $\theta_n$ is defined by
    \begin{equation*}
        (\beta_\omega R_n)^{1/\theta_n}=\beta_\omega R_n^{1/\theta}.
    \end{equation*}
    Such a choice is always possible by taking $R_n$ to be sufficiently small relative to $\beta_\omega$.
    That this indeed gives a lower bound for $\dimAs\theta E_\xi$ follows from \cite[Corollary~2.12]{zbl:1485.28006} (to be precise, this result is for the \emph{upper} Assouad spectrum, but the same proof works for the Assouad spectrum).
    The details of the version for the Assouad spectrum are also implicit in the proof of \cite[Theorem~2.1]{fhh+2019}).
    The remaining details are identical.
\end{proof}
As a quick application, the same argument allows us to extend \cref{p:uniform-qA} to the upper box dimension as well (alternatively, one may observe that the proof of \cite[Proposition~3.1]{fy2018b} is quantitative in the respective constants).
\begin{corollary}\label{c:uniform-box}
    Let $s=\sup_{\eta\in\Omega}\dimuB E_\eta$.
    Then for every $\epsilon>0$, there exists $C_\epsilon>0$ such that for all $0<R<1$ and $\eta\in\Omega$, $N_R(E_\eta)\leq C_\epsilon R^{-s-\epsilon}$.
\end{corollary}
\begin{remark}
    Since we can only guarantee that the fibre is large along a sequence of scales, we do not know if this result can be extended to the lower box dimension.
\end{remark}

\subsection{Assouad dimension and slices}\label{ss:adim}
We now turn our attention to the Assouad dimensions of dominated rectangular IFSs of affinities.
\begin{lemma}\label{l:fibre}
    Let $\{T_i\}_{i\in\mathcal{I}}$ be a dominated rectangular IFS of affinities with attractor $K$.
    Then
    \begin{equation}\label{e:fibre}
        \sup_{x\in\pi(K)}\dimA(\pi^{-1}(x)\cap K)\geq\max_{\eta\in\Omega}\bigl(\dimA E_\eta\bigr).
    \end{equation}
    Moreover, suppose the projected IFS $\{S_{\sigma}\}_{\sigma\in\mathcal{I}}$ satisfies the WSC.
    Then equality holds in \cref{e:fibre} and the supremum is attained as a maximum.
\end{lemma}
\begin{proof}
    As observed in \cref{ss:s-fibre}, $\{x_\eta\}\times E_\eta\subset K$ where $x_\eta=\lim_{n\to\infty}\eta_n(0)$.
    This directly gives \cref{e:fibre}.

    Now, assume that the projected IFS satisfies the WSC.
    Let $x\in K$ and set $s=\dimA(\pi^{-1}(x)\cap K)$.
    Applying \cref{p:wt-exist}, let $F\subset\R$ with $\mathcal{H}^s(F)>0$ be a weak tangent for $\pi^{-1}(x)\cap K$ with respect to the sequence of similarities $(U_j)_{j=1}^\infty$ with scaling ratios $\gamma_j$.
    Since the projected IFS $\{S_i\}_{i\in\mathcal{I}}$ satisfies the WSC,
    \begin{equation}\label{e:wsc-bound}
        \mathcal{M}_r\coloneqq\{S_\sigma:\sigma\in\Lambda_r,x\in S_\sigma(\pi(K))\}\qquad\text{satisfies}\qquad\sup_{r\in(0,1)}\#\mathcal{M}_r=M<\infty.
    \end{equation}

    Now for $t>0$, fix a packing $\{B(y_n,t)\}_{n=1}^{H_t}$ for $F$ with $y_n\in F$ and $H_t$ maximal.
    Since $\dimlB F>s-\epsilon$, for all $t$ sufficiently small, $H_t\geq (1/t)^{s-\epsilon}$.
    Now let $j$ be sufficiently large so that $d_{\mathcal{H}}(U_j(\pi^{-1}(x)\cap K)\cap B(0,1), F)\leq t/2$ and for all $\sigma\in\Lambda_{\gamma_j^{-1}}$, $\beta_\sigma\leq\gamma_j^{-1}t$.
    Such a choice is possible since $\beta_i<\alpha_i$ for all $i\in\mathcal{I}$.
    This implies that $U_j^{-1}(B(y_n,t))\cap \pi^{-1}(x)\cap K\neq\emptyset$ for all $n$.
    In particular, pigeonholing with respect to $\mathcal{M}_r$, there is some $y_t\in B(0,1)$ and $f_t\in\mathcal{M}_{\gamma_j^{-1}}$ so that for any $\eta\in\Omega$ with $\eta_k=f_t$ for some $k\in\N$,
    \begin{equation*}
        N_{\gamma_j^{-1} t}\bigl(E_\eta\cap B(y_t,\gamma_j^{-1}/2)\bigr)\geq\frac{H_{t}}{2M}\geq (1/t)^{s-2\epsilon}
    \end{equation*}
    by \cref{e:wsc-bound} for all $t$ sufficiently small.

    Finally, applying the argument to $\epsilon_n$ converging to zero, get a sequence $(t_n)_{n=1}^\infty$ converging to zero and let $\xi\in\Omega$ be chosen so that $\xi_{k_n}=f_{t_1}\circ\cdots\circ f_{t_n}$ for some $k_n\in\N$ and all $n\in\N$.
    Then the same argument from \cref{p:uniform-A} gives that $\dimA E_\xi\geq s-2\epsilon_n$ for all $n\in\N$, and therefore $\dimA E_\xi\geq s$.
    In particular, the supremum is attained as a maximum.
\end{proof}
Now, we prove our main result concerning Assouad dimensions.
When the projected IFS satisfies the WSC, applying \cref{l:fibre}, this also gives a formula for the Assouad dimension of a dominated rectangular IFS of affinities in terms of the Assouad dimension of the projection and the worst-case Assouad dimension of a fibre.
\begin{theorem}\label{t:dima}
    Let $\{T_i\}_{i\in\mathcal{I}}$ be a dominated rectangular IFS of affinities with attractor $K$.
    Then
    \begin{equation}\label{e:dimA-formula}
        \dimA K\geq \dimA \pi(K)+\max_{\eta\in\Omega}\bigl(\dimA E_\eta\bigr).
    \end{equation}
    Moreover, suppose the projected IFS $\{S_i\}_{i\in\mathcal{I}}$ satisfies the WSC.
    Then equality in \cref{e:dimA-formula} holds.
\end{theorem}
\begin{proof}
    Fix $\eta\in\Omega$: we first show that
    \begin{equation}\label{e:max-fibre}
        \dimA K\geq\dimA \pi(K)+\dimA E_\eta
    \end{equation}
    by constructing an appropriate weak tangent.

    Fix a weak tangent $\widehat{K}$ with positive $\mathcal{H}^{\dimA \pi(K)}$-measure for $\pi(K)$, and a weak tangent $\widehat{E}$ with positive $\mathcal{H}^{\dimA E_\eta}$-measure for the symbolic fibre $E_\eta$.
    For $n\in\N$, let $r_n$ denote the scaling ratio of $\eta_n$.

    Now let $\epsilon>0$ be arbitrary.
    Fix a similarity $U\colon\pi(K)\to\R$ with scaling ratio $\gamma$ so that
    \begin{equation*}
        d_{\mathcal{H}}(\widehat{K},U(\pi(K))\cap B(0,1))\leq\epsilon.
    \end{equation*}
    Next, let $n_0$ be sufficiently large so that $\beta_\sigma/\alpha_\sigma<\epsilon/\gamma$ for all $\sigma\in\mathcal{J}_{\eta_{n}}$ and all $n\geq n_0$.
    Then fix a similarity $V\colon E_\eta\to\R$ with scaling ratio $\kappa\geq \gamma\cdot r_{n_0}^{-1}$ so that
    \begin{equation*}
        d_{\mathcal{H}}(\widehat{E},V(E_\eta)\cap B(0,1))\leq\epsilon,
    \end{equation*}
    and let $n\geq n_0$ be such that $\gamma\cdot r_n^{-1}\beta_{\min}^{-1}\geq\kappa\geq\gamma\cdot r_n^{-1}$.
    In particular, since $\pi(K)$ is a self-similar set, as argued in \cref{l:wt-choice} there is a similarity $V_\epsilon$ with scaling ratio $\gamma\cdot r_n^{-1}$ and a similarity $h_\epsilon$ with scaling ratio in the interval $[\beta_{\min},1]$ so that
    \begin{equation*}
        d_{\mathcal{H}}(\widehat{E},V_\epsilon(E_\eta)\cap h_\epsilon(B(0,1)))\leq\epsilon.
    \end{equation*}
    Finally, $U_\epsilon=U\circ\eta_n^{-1}$ has scaling ratio $\gamma\cdot r_n^{-1}$.
    In particular, $W_\epsilon\colon K\to B(0,1)\times h_\epsilon(B(0,1))$ defined by $W_\epsilon(x,y)=(U\circ\eta_n^{-1}(x),V_\epsilon(y))$ is a similarity with scaling ratio $\gamma\cdot r_n^{-1}$.

    Thus taking an appropriate sequence $(\epsilon_n)_{n=1}^\infty$ tending to zero so that the scaling ratios of the $h_\epsilon$ converge, it follows that $\widehat{K}\times h(\widehat{E})$ is a weak pseudo-tangent for $K$, for some similarity map $h$.
    This gives \cref{e:max-fibre} by \cref{l:tan-lower}, and therefore \cref{e:dimA-formula}.

    We now prove the upper bound assuming the WSC.
    Let $\mathcal{Q}_R$ denote the collection of grid-aligned squares with side-length $R$.
    By the WSC assumption there is $M\in\N$ so that
    \begin{equation*}
        M=\sup_{R\in(0,1)}\sup_{Q\in\mathcal{Q}_R}\#\{S_\sigma:\alpha_\sigma\leq R<\alpha_{\sigma^-},S_\sigma(\pi(K))\cap \pi(Q)\neq\emptyset\}.
    \end{equation*}
    Let
    \begin{equation*}
        s_1=\dimA \pi(K)\text{ and }s_2=\max_{\eta\in\Omega}\bigl(\dimA E_\eta\bigr)
    \end{equation*}
    and let $\epsilon>0$ be arbitrary.
    (Note that $\pi(K)$ is Ahlfors--David regular since it satisfies the WSC, so in fact $s_1=\dimB\pi(K)$.)
    It suffices to show that $\dimA K\leq s_1+s_2+2\epsilon$.
    By \cref{p:uniform-A}, there is a constant $C>0$ so that for all $0<r\leq R<1$, $\eta\in\Omega$, and $x\in E_\eta$,
    \begin{equation}\label{e:uni}
        N_r(B(x,R)\cap E_\eta)\leq C\Bigl(\frac{R}{r}\Bigr)^{s_2+\epsilon}.
    \end{equation}

    Now fix $0<r\leq R<1$ and let $Q\in\mathcal{Q}_R$.
    Let $\gamma>1$ be chosen so that $\alpha_i^\gamma>\beta_i$ for all $i\in\mathcal{I}$.
    Let $\sigma\in\mathcal{I}^*$ be such that $\alpha_\sigma\leq R<\alpha_{\sigma^-}$ and $S_\sigma(\pi(K))\cap \pi(Q)\neq\emptyset$: we first cover
    \begin{equation*}
        Q\cap\bigcup_{\tau\in\mathcal{I}^*:S_\tau=S_\sigma}T_\tau(K).
    \end{equation*}
    Since $\dimA\pi(K)<s_1+\epsilon$, there is a constant $C_1>0$ and a family of balls $\{B_i\}_{i=1}^{\ell_1}\subset\R$ each with radius $r$ so that
    \begin{equation*}
        \pi(K\cap Q)\subset\bigcup_{i=1}^{\ell_1} B_i\quad\text{and}\quad \ell_1\leq C_1\Bigl({\frac{R}{r}}\Bigr)^{s_1+\epsilon}.
    \end{equation*}
    Thus by \cref{l:symb-fibre-limit} and the choice of $\gamma$, there is a constant $C>0$ so that
    \begin{equation}\label{e:ph-close}
        p_{\mathcal{H}}\left(\bigcup_{\tau\in\mathcal{I}^*:S_\tau=S_\sigma}T_\tau(K),E_\eta\right)\leq C\max\{\beta_\tau:\tau\in\mathcal{I}^*:S_\tau=S_\sigma\}\leq C\alpha_\sigma^\gamma\leq CR^\gamma.
    \end{equation}
    where $\eta\in\Omega$ is any sequence with $\eta_k = S_\sigma$ for some $k\in\N$.
    Moreover, it suffices in the definition of the Assouad dimension to consider scales $r\geq CR^{\gamma}$: by \cite[Proposition 3.7]{fy2018b}, for any $\theta\in(0,1)$, there is a $\theta'>1/\gamma$ such that $\dimAs{\theta'} K\geq\dimAs\theta K$.
    (Alternatively, one can directly apply the covering argument from \cref{t:dimas}).
    Thus by \cref{e:uni} combined with \cref{e:ph-close}, there is a constant $C_2>0$ and family of balls $\{\widehat{B}_i\}_{i=1}^{\ell_2}\subset\R$ each with radius $r$ so that
    \begin{equation*}
        Q\cap\bigcup_{\tau\in\mathcal{I}^*:S_\tau=S_\sigma}T_\tau(K)\subset\bigcup_{i=1}^{\ell_1}\bigcup_{j=1}^{\ell_2}2B_i\times\widehat{B}_i\quad\text{and}\quad\ell_2\leq C_2\Bigl(\frac{R}{r}\Bigr)^{s_2+\epsilon}
    \end{equation*}
    where $2B_i$ is the ball with the same centre as $B_i$ and double the radius.
    Note that $C_1$ and $C_2$ depend only on the governing IFS.
    Thus
    \begin{equation*}
        N_r(K\cap Q)\leq 3MC_1C_2\Bigl(\frac{R}{r}\Bigr)^{s_1+s_2+2\epsilon}
    \end{equation*}
    where $3$ is the doubling constant in $\R$, and the result follows.
\end{proof}
\begin{remark}\label{r:awsc-req}
    In general, the WSC assumption here is required.
    For example, consider the IFS given by the maps
    \begin{equation*}
        T_1(x,y)=(\beta x,\alpha y)\text{ and }T_2(x,y)=(\beta x+(1-\beta),\alpha y+(1-\alpha))
    \end{equation*}
    where $0<\alpha\leq 1/2<\beta<1$.
    Denote the attractor of such a set by $K$.
    Such carpets are often referred to as Przytycki--Urbański sets \cite{pu1989}.

    In \cite[Section~2.1]{fj2017}, the dimension of such carpets are computed for various values of $\alpha$ and $\beta$.
    In particular, if $1/\beta$ is a Garsia number, i.e.~a real algebraic integer with norm $2$ and Galois conjugates lying strictly outside the complex unit disc, then
    \begin{equation*}
        \dimA K=1+\frac{\log 2\beta}{-\log \alpha}.
    \end{equation*}
    However, it is easy to check that the IFS $\{x\mapsto \beta x,x\mapsto \beta x+(1-\beta)\}$ has no exact overlaps.
    Thus the lower bound from \cref{e:dimA-formula} gives the value $1$.
\end{remark}
\subsection{Quasi-Assouad dimension}
We now turn our attention to the Assouad spectrum, for values of $\theta$ sufficiently close to $1$.
In general, for smaller values of $\theta$, determining a precise formula for the Assouad spectrum seems to be rather complicated since the Assouad spectrum is highly sensitive to the logarithmic eccentricity ratios $(\log\alpha_i)/(\log\beta_i)$.
\begin{theorem}\label{t:dimas}
    Let $\{T_i\}_{i\in\mathcal{I}}$ be a dominated rectangular IFS of affinities with attractor $K$.
    Then for all $\max_{i\in\mathcal{I}}\frac{\log\alpha_i}{\log\beta_i}\leq\theta<1$,
    \begin{equation}\label{e:dimAs-formula}
        \dimAs\theta K\geq\dimH \pi(K)+\max_{\eta\in\Omega}\bigl(\dimAs\theta E_\eta\bigr)
    \end{equation}
    Moreover, if the projected IFS $\{S_i\}_{i\in\mathcal{I}}$ satisfies the AWSC, then equality in \cref{e:dimAs-formula} holds.
\end{theorem}
\begin{proof}
    Let $\theta_0=\max_{i\in\mathcal{I}}\frac{\log\alpha_i}{\log\beta_i}$ and let $\theta\in(\theta_0,1)$ (the result follows at $\theta_0$ by continuity of the Assouad spectrum).
    This is equivalent to requiring that $\alpha_\sigma^{1/\theta}\geq\beta_\sigma$ for all $\sigma\in\mathcal{I}^*$.
    Let $\epsilon>0$: we will first show that, assuming the AWSC,
    \begin{equation*}
        \dimAs\theta K\leq\dimuB\pi(K)+\max_{\eta\in\Omega}\bigl(\dimAs\theta E_\eta\bigr)+\epsilon,
    \end{equation*}
    recalling that $\pi(K)$ is a self-similar set and therefore has $\dimH\pi(K)=\dimuB\pi(K)$.

    Let $r>0$ be arbitrary and consider a grid-aligned square $Q$ with side-length $r$.
    By the AWSC assumption, for all $x\in\pi(K)$ and $r$ sufficiently small, with
    \begin{equation*}
        \mathcal{M}(x,r)\coloneqq\{S_\sigma:\sigma\in\Lambda_r,S_\sigma(\pi(K))\cap B(x,r)\neq\emptyset\}
    \end{equation*}
    we have $\sup_{x\in K}\#\mathcal{M}(x,r)\leq \bigl(r^{1-1/\theta}\bigr)^{\epsilon}$.
    Moreover, for fixed $S_\sigma\in\mathcal{M}(x,r)$, by \cref{p:uniform-qA} for all $r$ sufficiently small (depending only on the IFS) there is a family of balls $\{B_i\}_{i=1}^N$ each with radius $r^{1/\theta}$ such that for any $\eta\in\Omega$ with $\eta_n=S_\sigma$, $\{B_i\}_{i=1}^N$ is a cover for $\overline{\pi}(Q)\cap E_\eta$, and
    \begin{equation*}
        N\leq\Bigl(\frac{r}{r^{1/\theta}}\Bigr)^{\dimAs\theta E_\eta+\epsilon}.
    \end{equation*}
    In particular, since $\theta>\theta_0$, by \cref{l:symb-fibre-limit} for all $r$ sufficiently small
    \begin{equation}\label{e:close}
        p_{\mathcal{H}}\left(\bigcup_{\tau\in\mathcal{I}^*:S_\tau=S_\sigma}T_\tau(K),E_\eta\right)\leq r^{1/\theta}.
    \end{equation}
    Similarly, let $\{U_i\}_{i=1}^M$ be a cover for $S_\sigma(\pi(K))$ where each $U_i$ is a ball with radius $r^{1/\theta}$ so that
    \begin{equation*}
        M\leq\Bigl(\frac{r}{r^{1/\theta}}\Bigr)^{\dimuB\pi(K)+\epsilon}.
    \end{equation*}
    Such a cover exists by taking the image under $S_\sigma$ of a cover for $\pi(K)$ at scale $r^{1/\theta-1}$.
    Therefore by \cref{e:close},
    \begin{equation*}
        Q\cap\bigcup_{\substack{\tau\in\mathcal{I}^*\\S_\tau=S_\sigma}}T_\tau([0,1]^2)\subset\bigcup_{i=1}^N\bigcup_{j=1}^M 2B_i\times U_j.
    \end{equation*}
    Thus taking a union over $\mathcal{M}(x,r)$,
    \begin{equation}\label{e:upper}
        N_{r^{1/\theta}}(K\cap Q)\leq 3\Bigl(\frac{r}{r^{1/\theta}}\Bigr)^{\epsilon}\Bigl(\frac{r}{r^{1/\theta}}\Bigr)^{\dimAs\theta E_\eta+\epsilon}\Bigl(\frac{r}{r^{1/\theta}}\Bigr)^{\dimuB\pi(K)+\epsilon}
    \end{equation}
    where $3$ is the doubling constant in $\R$, from which the upper bound follows.

    To obtain the converse inequality, write $\alpha=\sup_{\eta\in\Omega}\dimAs\theta E_\eta$, let $\epsilon>0$, and let $\gamma\in\Omega$ be such that
    \begin{equation*}
        \dimAs\theta E_{\gamma}\geq\alpha-\epsilon/2.
    \end{equation*}
    Then get $x\in E_\gamma$ and $r$ arbitrarily small so that
    \begin{equation*}
        N_{r^{1/\theta}}\bigl(B(x,r)\cap E_\gamma\bigr)\geq\left(\frac{r}{r^{1/\theta}}\right)^{\alpha-\epsilon}.
    \end{equation*}
    Finally, let $k\in\N$ and $\tau\in\mathcal{I}^k$ be such that $\gamma_k=S_\tau$ and $r\alpha_{\min}\leq\alpha_\sigma\leq r$.
    Since $\theta>\theta_0$, if $r$ is chosen to be sufficiently small, $\beta_\sigma\leq r^{1/\theta}$ for all $\sigma\in\mathcal{J}_{S_\tau}$.
    Thus for $r$ sufficiently small,
    \begin{equation*}
        N_{r^{1/\theta}}\bigl((B(x,r)\times B(S_\sigma(x_0),r))\cap K\bigr)\geq\left(\frac{r^{1/\theta}}{r}\right)^{\alpha-\epsilon}\left(\frac{r^{1/\theta}}{r}\right)^{\dimH \pi(K)-\epsilon}
    \end{equation*}
    for some fixed $x_0\in\pi(K)$.
    But $\epsilon>0$ and $r>0$ can be chosen to be arbitrarily small, so that $\dimAs\theta K\geq\alpha+\dimH\pi(K)$, as claimed.
\end{proof}
\begin{remark}
    The proof of the lower bound for \cref{e:dimAs-formula} is simpler than the Assouad dimension lower bound in \cref{t:dima} since the quasi-Assouad dimension only witnesses the box dimension of $\pi(K)$, whereas the Assouad dimension picks up the Assouad dimension of $\pi(K)$.
    The box dimension of $\pi(K)$ always exists, but we must carefully move tangents for $\dimA\pi(K)$ in \cref{t:dima}.
    This is similar to the phenomenon in which a self-similar set $K$ in $\R$ which satisfies the AWSC but not the WSC has $\dimH K=\dimqA K$ but $\dimA K=1$.
\end{remark}
In fact, since any self-similar set in $\R$ is the attractor of a dominated rectangular IFS of affinities, we obtain the following corollary:
\begin{corollary}\label{c:ss-qa}
    Let $\{S_i\}_{i\in\mathcal{I}}$ be an IFS of similarities in $\R$ satisfying the AWSC with attractor $K$.
    Then $\dimqA K=\dimH K$.
\end{corollary}
This was previously observed under the ESC (see \cite[Theorem~7.3.1]{fra2021} and the discussion which follows it), which is a special case of our result by \cref{p:esc-class}.

Combining \cref{c:ss-qa} with \cref{t:dimas}, we obtain our main result on quasi-Assouad dimensions.
\begin{corollary}\label{c:dimqa}
    Let $\{T_i\}_{i\in\mathcal{I}}$ be a dominated rectangular IFS of affinities with attractor $K$.
    Suppose the projected IFS $\{S_i\}_{i\in\mathcal{I}}$ satisfies the AWSC.
    Then
    \begin{equation*}
        \dimqA K=\dimqA \pi(K)+\max_{\eta\in\Omega}\bigl(\dimqA E_\eta\bigr).
    \end{equation*}
\end{corollary}
\subsection{Distinct Assouad and quasi-Assouad dimensions}\label{ss:ques}
Our original motivation for studying this family of examples was in the context of \cite[Question~17.5.4]{fra2021}, which asks whether or not it can happen that $\dimqA K<\dimA K$ for a self-affine set $K$ satisfying the strong separation condition.

We show that $\dimqA K<\dimA K$ is possible: to be precise, for any $\epsilon>0$, we construct a self-affine set $K$ satisfying the strong separation condition such that $\dimqA K<\epsilon$ and $\dimA K\geq 1$.
First, we observe the following application of \cref{c:dimqa}:
\begin{corollary}\label{c:dim-no-fibre}
    Let $\{T_i\}_{i\in\mathcal{I}}$ be a dominated rectangular IFS of affinities such that the projected IFS $\{S_i\}_{i\in\mathcal{I}}$ satisfies the ESC with similarity dimension $s\leq 1$.
    Then $\dimH K=\dimqA K=\dimqA \pi(K)$.
\end{corollary}
\begin{proof}
    Since the projected IFS has no exact overlaps, $\dimqA E_\eta=0$ for all $\eta\in\Omega$.
    Since $s\leq 1$, the projected IFS also satisfies the AWSC by \cref{p:esc-class}.
    Then apply \cref{c:dimqa}.
\end{proof}
In general, the ESC holds for typical parameters for parametrized families of self-similar sets (see, for example, \cite[Theorem~1.7]{hoc2014}), whereas this is not true for the OSC.
In particular, it is well-known that there exist examples which satisfy the ESC with similarity dimension arbitrarily close to $0$, but not the OSC.

Here, we will construct an explicit example of this phenomenon for the convenience of the reader.
Let $\epsilon\in(0,1)$ be arbitrary and fix $N\in\N$ with $N\geq 4$ so that $(\log 3)/(\log N)\leq\epsilon$.
First, for $t\in(0,1/N)$, consider the IFS $\Phi_t=\{S_0,S_t,S_{1-1/N}\}$ where $S_j(x)=\frac{x}{N}+j$ for $j\in\{0,t,1-1/N\}$.
Observe that if $t\notin\Q$, then the IFS has no exact overlaps.
Moreover, the same argument as \cite[Theorem~1.6]{hoc2014} shows that the IFS $\Phi_t$ satisfies the ESC for all $t\notin\Q$.

Suppose $\Phi_t$ satisfies the OSC, i.e.~there is an open set $U\subset(0,1)$ such that
\begin{equation*}
    U\supset S_1(U)\cup S_t(U)\cup S_{1-1/N}(U)
\end{equation*}
disjointly.
Then observe that if $\delta>0$ is a constant such that there is some $x$ with $B(x,\delta/2)\subset U$, then
\begin{equation}\label{e:OSC-bd}
    |S_\sigma(0)-S_\tau(0)|\geq\frac{\delta}{N^n}
\end{equation}
for all $\sigma\neq\tau\in\mathcal{I}^n$.

Let $C$ denote the attractor of the IFS $\{S_0,S_{1-1/N}\}$, i.e.~$C$ is the usual Cantor set with subdivision ratios $1/N$.
Now let $t_0$ be an irrational number which is approximated well by left endpoints of the level $n$ intervals of the Cantor set, i.e.~for infinitely many $n\in\N$, there are $d_1,\ldots,d_n\in\{0,N-1\}$ so that
\begin{equation*}
    \Bigl\lvert t_0-\frac{\sum_{i=1}^n d_i N^{i-1}}{N^n}\Bigr\rvert\leq\frac{\delta_n}{N^n}
\end{equation*}
with the $\delta_n$ converging to $0$.
For example, let $t_0\in(0,1/N)\setminus\Q$ be a number with base $N$ representation consisting of arbitrarily long sequences of $0$s or $(N-1)$s.
Since $t_0$ is not an endpoint of a level $n$ interval for all $n\geq 1$, it follows immediately from \cref{e:OSC-bd} that $\Phi_{t_0}$ cannot satisfy the OSC.

Then the attractor $K$ of the IFS
\begin{align*}
    T_1(x,y)&=\Bigl(\frac{x}{N},\frac{y}{N+1}\Bigr)\\
    T_2(x,y)&=\Bigl(\frac{x}{N}+t_0,\frac{y}{N+1}+\frac{N-1}{2(N+1)}\Bigr)\\
    T_3(x,y)&=\Bigl(\frac{x}{N}+1-\frac{1}{N},\frac{y}{N+1}+\frac{N}{N+1}\Bigr)
\end{align*}
satisfies the strong separation condition and has, for $N$ sufficiently large,
\begin{equation*}
    \dimqA K=\dimH K=\frac{\log 3}{\log N}\leq\epsilon<1\leq\dimA K.
\end{equation*}
by \cref{c:dim-no-fibre}.
We recall that $\dimA K\geq\dimA\pi(K)$ since a self-similar image of $\pi(K)$ is always contained in some weak tangent of $K$.
See the discussion in the introduction following \cref{ic:Assouad-formula}.
\begin{remark}\label{r:it-wsc-fail}
    Similar arguments also show that there are arbitrarily small parameters $0<\lambda_1,\lambda_2\leq\rho\coloneqq(3-\sqrt{5})/2$ such that the IFS $\Phi$ defined in \cref{e:lambda-ifs} fails the WSC.
    To see this, let $0<r<\rho^2$ be small, fix $\lambda_1\in(r/\rho,\sqrt{r})$, and let $\lambda_2=r/\lambda_1$.
    Observe that
    \begin{equation*}
        \Phi\supset\{x\mapsto rx, x\mapsto rx+(1-r), x\mapsto rx+\lambda_1-\lambda_1^2, x\mapsto rx+\lambda_1-r\}.
    \end{equation*}
    In particular, if $r-\lambda_1^2$ is close to a left endpoint $b$ of a level $n$ interval of the IFS $\{x\mapsto rx, x\mapsto rx+1-r\}$, then $\lambda_1-\lambda_1^2$ is close to $rb+\lambda_1-r$, which is the left endpoint of a level $n+1$ interval of the IFS $\Phi$.
    Moreover, it is well-known that if the WSC holds, then there is a constant $\delta>0$ such that for all $f_1,\ldots,f_n$ and $g_1,\ldots,g_n$ in $\Phi$, either $f_1\circ\cdots\circ f_n(0)=g_1\circ\cdots\circ g_n(0)$ or
    \begin{equation}\label{e:WSC-bd}
        |f_1\circ\cdots\circ f_n(0)-g_1\circ\cdots\circ g_n(0)|\geq\frac{\delta}{r^n}.
    \end{equation}
    (see, for example, \cite[Theorem~1]{zer1996}).
    In particular, one may choose $\lambda_1$ arbitrarily close to $\sqrt{r}$ so that \cref{e:WSC-bd} fails for any $\delta>0$ for the IFS $\Phi$.
    This choice of $\lambda_1$ also guarantees that $\lambda_2$ is arbitrarily close to $\sqrt{r}$, and hence is arbitrarily small as well.
\end{remark}
\begin{acknowledgements}
    JMF was supported by Leverhulme Trust Research Project Grant (RPG-2019-034) and RSE Sabbatical Research Grant (70249).
    AR was supported by EPSRC Grant EP/V520123/1 and the Natural Sciences and Engineering Research Council of Canada.
    The authors thank Antti Käenmäki and Sascha Troscheit for various discussions involving Assouad dimensions and related topics.
    They also thank Aleksi Pyörälä, Adam Śpiewak, and the anonymous referee for (independently) pointing out an error in the definition of the symbolic slice.
\end{acknowledgements}
\end{document}